\documentclass[12pt, oneside, leqno]{article}
\usepackage[T1]{fontenc}
\usepackage{amsmath,amsthm}
\usepackage{amssymb,latexsym}
\usepackage[a4paper,scale={0.7,0.8}]{geometry}
\usepackage{enumitem}
\newtheorem{thm}{Theorem}
\newtheorem{cor}{Corollary}
\newtheorem{prop}{Proposition}
\newtheorem{lem}{Lemma}
\newtheorem{prob}{Problem}
\numberwithin{equation}{section}
\theoremstyle{definition}

\newtheorem{rem}{Remark}

\setenumerate{label={(\arabic*)},nolistsep}

\begin{document}
\baselineskip=16pt
\author{Piotr  Nowak\\
Mathematical Institute,  University of Wroc{\l}aw\\
Pl. Grunwaldzki 2/4, 50-384 Wroc{\l}aw, Poland\\
E-mail: nowak@math.uni.wroc.pl
}
\title{Some properties of the moment estimator of shape parameter
for the gamma distribution}
\date{}
\maketitle
%%%%%%%%%%%%%%%%%%%%%%%%%%%%%%%%%%%%%%%%%%%%%%%%%%%%%%%%%
\begin{abstract}
Exact distribution of the moment estimator of shape parameter for the 
gamma distribution for small samples is derived. Order preserving properties
of this estimator are presented.

\end{abstract}
%%%%%%%%%%%%%%%%%%%%%%%%%%%%%%%%%%%%%%%%%%%%%%%%%%%%%%%%%
\section{Introduction and preliminaries}

It is well known that the gamma distribution has wide 
application, for example in meteorology to describe the 
distribution of rainfall.

The moment estimators  $\hat\alpha$ and $\hat\lambda$
of the parameters $\alpha$ and $\lambda$ of the 
gamma distribution with density

$$
f(x;\alpha,\lambda) = \frac{1}{\Gamma(\alpha) \lambda^\alpha} x^{\alpha -1} \exp (-x/\lambda), \ \ x>0, \ \alpha>0, \lambda>0
$$
for a random sample $X_1,\dots,X_n$ are
 
\begin{equation}
\label{mm_est}
\hat\alpha=\frac{\bar X^2}{S^2} \quad
\textrm{and} \quad
\hat\lambda=\frac{S^2}{\bar X},
\end{equation}
where $\bar X$ and $S^2$ 
are  sample mean and sample biased variance respectively.

In this paper we are mainly interested in order preserving property of
the estimator $\hat\alpha$.
In general, we say that the estimator $\hat\theta$  based
on the sample $X_1,\dots,X_n$ from population with density 
$f(\cdot;\theta)$ is increasing in $\theta$ with respect
to the order $\prec$  if $\hat\theta_1\prec\hat\theta_2$,
whenever $\theta_1<\theta_2$, where $\hat\theta_1$
($\hat\theta_2$) is the estimator based on sample 
from density $f(\cdot;\theta_1)$ ($f(\cdot;\theta_2)$).
We shall consider when $\prec$ is one of the following
orders: $\leq_{\textup{st}}$, $\leq_{\textup{disp}}$, 
$\leq_{\textup{*}}$. For details and definitions of these stochastic orders 
we refer the reader to Shaked and Shanthikumar~\cite{shaked_st:or}.

Nowak \cite{nowak:2} proved the following theorem.

\begin{thm}\label{thm_alpha}
The estimators $\hat\alpha$ and $\hat\lambda$ are stochastically increasing
respectively in $\alpha$ and~$\lambda$.
\end{thm}

Deriving exact distributions of  moment estimators $\hat \alpha$
and $\hat \lambda$ is very tedious for $n>2$.
Now we find the distribution of $\hat\alpha$ for $n=2$. First we
prove the more general fact.

\begin{lem}
Suppose that $X_1\sim G(\lambda,\alpha_1)$ and 
$X_2\sim G(\lambda,\alpha_2)$, $X_1$ and $X_2$ 
are independent random variables.
Then $Z=\dfrac{S^2}{\bar X^2}$ has density of the form
\begin{align*}
&h(z;\alpha_1,\alpha_2)=
1/B(\alpha_1,\alpha_2)
2^{-(\alpha_1+\alpha_2)} 
z^{-1/2}\times \\& \times
\left[(1-\sqrt{z})^{\alpha_1-1}(1+\sqrt{z})^{\alpha_2-1}+
(1+\sqrt{z})^{\alpha_1-1}(1-\sqrt{z})^{\alpha_2-1}\right].
\end{align*}
\end{lem}
\begin{cor}
Let $X_1, X_2\sim G(\lambda,\alpha)$, $Y\sim Ex(1)$, 
$X_1$, $X_2$ and $Z$ are independent random variables.
Then
\begin{equation}
\label{exp_char}
\frac{X_1^2+X_2^2}{(X_1+X_2)^2}=_{\textup st} 
\frac{(X_1+Y)^2+X_2^2}{(X_1+Y+X_2)^2}.
\end{equation}
\end{cor}
\begin{prob}
Equation (\ref{exp_char}) yields some interesting property of the exponential
distribution. There arises a question: Does this equation
characterize the exponential distribution?
\end{prob}

\begin{cor}
For $n=2$ the statistics
$1/\hat\alpha$ has the beta distribution  $B(1/2,\alpha)$.
Hence $\hat\alpha$ has monotone likelihood ratio.
\end{cor}

Now we examine dispersive ordering of the $\hat\alpha$ estimator.
First we notice that $1/\hat\alpha$ is not dispersive monotone in $\alpha$.
We conclude it from Theorem 3.B.14 of
Shaked and Shanthikumar \cite{shaked_st:or} since  the support of the
estimator $1/\hat\alpha$ is the finite interval $(0,n-1)$.

\begin{thm}
The estimator $1/\hat\alpha$ is not dispersively monotone in $\alpha$.
\end{thm}

We now consider dispersive ordering of the estimator $\hat\alpha$.
From previous theorem it does not follow that $\hat\alpha$ is not
dispersive monotone in $\alpha$, though the function $1/x$
is decreasing and convex for $x>0$.  We should notice that  
Theorem 3.B.10(b) of  Shaked and Shanthikumar \cite{shaked_st:or}
is not valid here. 
For example, in the following theorem we prove that for $n=2$ the estimator
$\hat\alpha$ is dispersively increasing in $\alpha$.

\begin{thm}
For $n=2$ the estimator $\hat\alpha$ is dispersively increasing in $\alpha$.
\end{thm}

\begin{proof}
Let $f^{\alpha}$ denotes density of the estimator $\hat\alpha$.
We know that $1/\hat\alpha$ has the beta distribution 
$\mathcal{B}(1/2,\alpha)$, so 
$$
f^{\alpha}(x)=\frac{1}{B(1/2,\alpha)}x^{-3/2}(1-1/x)^{\alpha-1}, \ x>1.
$$
In order to prove that
$\hat\alpha_1\leq_{\textrm{disp}}\hat\alpha_2$
whenever $\alpha_1<\alpha_2$ it suffices to prove that 
$$S^{-}(f^{\alpha_1}_{c}(x+1)-f^{\alpha_2}(x+1))
\leq 2 \quad \textrm{for\ every} \ c>0$$
with the sign sequence being $-,+,-$ in the case of equality
and $\hat\alpha_1\leq_{\textrm{st}}\hat\alpha_2$ (see Theorem~2.6 
and Remark 2.1 of Shaked \cite{shaked}).\\
Fix $c>0$. When $x>c$ then
\begin{flalign*}
&f^{\alpha_1}_{c}(x+1)-f^{\alpha_2}(x+1)  
=\frac{1}{B(1/2,\alpha_2)}(x+1)^{-3/2}(1-1/(x+1))^{\alpha_2-1}\times \\&
\times\big[A (x+1)^{3/2}(x-c+1)^{-3/2}
(1-1/(x-c+1))^{\alpha_1-1}/(1-1/(x+1))^{\alpha_2-1}-1\big],
\end{flalign*}
where $A=B(1/2,\alpha_2)/B(1/2,\alpha_1)$. Thus
$
S^{-}(f^{\alpha_1}_{c}-f^{\alpha_2})= S^{-}(h(x)),
$ where 
\begin{flalign*}
h(x)&=\log A + 3/2\log(x+1) - 3/2\log(x-c+1)+\\&
+(\alpha_1-1)\log(1-1/(x-c+1))
 - (\alpha_2-1)\log (1-1/(x+1)).
\end{flalign*}
Now we show that $h$ has exactly one maximum when $\alpha_1>1$ and 
the sign sequence is $-,+,-$. In the other case $h$ is decreasing.
After calculating we have that $h'(x)=0$ if and only if 
$
-w(x)/(2x(x+1)(x-c)(x-c+1))=0,
$ where 
$$
w(x)=(2\alpha_2-2\alpha_1+3c)x^2+(2\alpha_2-2\alpha_1+4c-4\alpha_2c-3c^2)x+
2c-2\alpha_2c-2c^2+2\alpha_2c^2.
$$ So we must show that $w$ has only one root in the interval $(c,\infty)$ 
if $\alpha_1>1$ and no roots when $\alpha_{1}\in(0,1]$.
Define 
$\bar{w}(x)=w(x+c)$. Therefore
$$\bar{w}(x)=(2(\alpha_2-\alpha_2)+3c)x^2+(2(\alpha_2-\alpha_1)+4c(1-\alpha_1)+3c^2)x
+2(1-\alpha_1)(c+c^2).$$
It is easy to see that for $\alpha_1>1$, $\bar{w}$ has two different roots $x_1$
and $x_2$ that $x_1 x_2 <0$ and for $\alpha_1\in(0,1]$ $\bar{w}$ has no 
roots in $(0,\infty)$. At the end we must  show that 
 the sign sequence is $-,+,-$ when $S^{-}(h)=2$. For $\alpha_1>1$
$$
\lim_{x\to c^+} h(x) = -\infty, \ \lim_{x\to \infty} h(x) = \log A <0.
$$
Combining these above facts with Theorem \ref{thm_alpha} we end the proof.
\end{proof}
\begin{rem}
The above theorem  does not hold in general  for $n\geq 3$, but
direct calculating is impossible due to occurrence of
hyper elliptic integrals. For example, numerical 
calculations show, that for $n=3$, $\hat\alpha_1\not<_{\textrm{disp}}\hat\alpha_2$
if $\alpha_1=1/5$ and $\alpha_2=1/4$. Then $G^{-1}(x)-F^{-1}(x)$
has the local maximum at $x\approx 0.72$ and the local minimum at $x\approx 0.85$.
\end{rem}
%%%%%%%%%%%%%%%%%%%%%%%%%%%%%%%%%%%%%%%%%%%%%%%%%%%%%%%%%%%%%%%%%%%%%
%%%%%%%%%%%%%%%%%%%%%%%%%%%%%%%%%%%%%%%%%%%%%%%%%%%%%%%%%%%%%%%%%%%%%
%%%%%%%%%%%%%%%%%%%%%%%%%%%%%%%%%%%%%%%%%%%%%%%%%%%%%%%%%%%%%%%%%%%%%
%%%%%%%%%%%%%%%%%%%%%%%%%%%%%%%%%%%%%%%%%%%%%%%%%%%%%%%%%%%%%%%%%%%%%
Deriving exact distribution when $n=3$ is more tedious.  We should 
add, that the estimator $\hat\alpha$ is closely related to the student 
statistic $\sqrt{n}\bar X/S$. The properties of this statistic
was very intensively studied in the literature, especially 
under non-normal conditions. First,
we find joint distribution of the vector $(\bar X, S)$.
One can proof the following lemma, see for example Craig \cite{craig}.
\begin{lem}
If $f$  is  a density on $(0,\infty)$, then the joint distribution 
$(\bar X, S)$  for $n=3$ is given by formula
$$
F(\bar x, s)=\left \{
\begin{array}{llr}
18s&\int_{\bar x-s\sqrt{2}}^{\bar x+s\sqrt{2}} \frac{1}{R} f(x_1)
f(\frac{3\bar x -x_1+R}{2})f(\frac{3\bar x -x_1-R}{2})dx_1, &
0\leq s\leq \bar x\sqrt{2}/2,\\[2ex]
18s &\left[\int_{0}^{\frac{3\bar x -\sqrt{6s^2-3\bar x^2}}{2}}+
\int^{\bar x+s\sqrt{2}}_{\frac{3\bar x +\sqrt{6s^2-3\bar x^2}}{2}}\right]&\\[2ex]
&\frac{1}{R} f(x_1) f(\frac{3\bar x -x_1+R}{2})f(\frac{3\bar x -x_1-R}{2})dx_1, &
 \bar x\sqrt{2}/2 \leq s \leq \bar x \sqrt{2},
\end{array}
\right.
$$
\nopagebreak[4]
where $R=\sqrt{6s^2-3(x_1-\bar x)^2}$.
\nopagebreak[4]
\end{lem}
If we have the joint distribution $(\bar X, S)$ it is easy to 
derive the  distribution of the statistics $T=\bar X/S=\sqrt{\hat \alpha}$
from the formula 
$$f_T(t)=\int u F(ut,u)du.$$
\begin{prop}
For $n=3$
the  cumulative distribution function of $T$
 for the exponential distribution with
density $f(x)=e^{-x}$, $x>0$ 
is given by
$$
F(t)=\Bigg\{\begin{array}{lr}
1-\frac{2 \pi }{3 \sqrt{3} t^2},& t\geq \sqrt{2},\\
-\frac{\sqrt{2-t^2}}{\sqrt{3} t}+\frac{\pi }{3 \sqrt{3} t^2}-\frac{2 \arcsin\left(\frac{t}{\sqrt{2}}\right)}{\sqrt{3} t^2}+1, &
\frac{\sqrt{2}}{2}\leq t\leq \sqrt{2}.  
\end{array}
$$ 
\end{prop}
On the other hand it is not easy to derive in the general the
distribution of~$T$. For example, if $\alpha\in(0,1)$ above 
integrals can be calculated only by numerical methods. For $n=3$,
after a bit algebra we can also derive the distribution of $T$ 
 for the  gamma  distribution with shape parameter $\alpha=2$.

\begin{prop}
For $n=3$ the  cumulative distribution function of $T$ 
for the gamma distribution with
density $g(x)=xe^{-x}$, $x>0$ is given by
$$
G(t)=\left\{\begin{array}{lr}
1-\frac{10 \pi  \left(4 t^2-3\right)}{27 \sqrt{3} t^4},& t\geq \sqrt{2},\\[2ex]
\frac{\sqrt{2-t^2} \left(-33 t^4-13 t^2+8\right) +5 \pi  t \left(4 t^2-3\right)-30 t \left(4 t^2-3\right) \arcsin\left(\frac{t}{\sqrt{2}}\right)}{27 \sqrt{3} t^5}+1
, &
\frac{\sqrt{2}}{2}\leq t\leq \sqrt{2}.  
\end{array}\right.
$$ 
\end{prop}

From previous propositions it follows, that the statistic $T$
is not monotone with respect to the star order $\leq_{*}$, 
that is the function $G^{-1}F(x)/x$ is not monotone in~$x>0$. To~see it,
let us  calculate:\\[2ex]
$
1) \
\frac{G^{-1}F(x)}{x}\bigg|_{x=\frac{11}{10}}=
\frac{10}{33}
\sqrt{\frac{11 \left(220 \pi +\sqrt{10 \pi  \left(-891 \sqrt{79}-560 \pi +16200 \arccos\left(\frac{11}{10 \sqrt{2}}\right)\right)}\right)}{33 \sqrt{79}+200 \pi -600 \arccos
   \left(\frac{11}{10 \sqrt{2}}\right)}}
 \approx 1.32686\\[2ex]
$
\noindent
$
2) \
\frac{G^{-1}F(x)}{x}\bigg|_{x=\frac{\sqrt{6}}{2}}=
\sqrt{\frac{40 \pi +2 \sqrt{10 \pi  \left(22 \pi -27 \sqrt{3}\right)}}{27 \sqrt{3}+18 \pi }}\approx 1.31502\\[2ex]
$
\noindent
$
3) \
\frac{G^{-1}F(x)}{x}\bigg|_{x=\sqrt{2}}=
\frac{\sqrt{40+2\sqrt{130}}}{6} \approx 1.32081\\
$

Since the star order is preserved under increasing function,
we have the following corollary.
\begin{cor}
For $n=3$ the estimator $\hat\alpha$  is not increasing in $\alpha$ with respect to the star order.
\end{cor}

The star order is closely related to the  dispersive order, since
$$
X\leq_{\textup{*}} Y \quad
\textup{iff} \quad  \log X\leq_{\textup{disp}} \log Y.
$$
Using the same technique as in Theorem 3 we can proof the following
proposition.

\begin{prop}
For $n=2$ the estimator $\hat\alpha$  is increasing in $\alpha$ with respect to the star order.
\end{prop}

Applying the property that for nonnegative random variables
such $X\leq_{\textup{st}} Y$ and
$X\leq_{\textup{*}}~Y$ implies the ordering $X\leq_{\textup{disp}}~Y$ (see, for example, Ahmed et al.~\cite{bart}) we can  also deduce from Theorem 1 and
Proposition 3 that $\hat\alpha$
is  dispersively monotone for $n=2$.
%%%%%%%%%%%%%%%%%%%%%%%%%%%%%%%%%%%%%%%%%%%%%%%%%%%%%%%%%%%%%%%%%%%%%
%%%%%%%%%%%%%%%%%%%%%%%%%%%%%%%%%%%%%%%%%%%%%%%%%%%%%%%%%%%%%%%%%%%%%
%%%%%%%%%%%%%%%%%%%%%%%%%%%%%%%%%%%%%%%%%%%%%%%%%%%%%%%%%%%%%%%
%%%%%%%%%%%%%%%%%%%%%%%%%%%%%%%%%%%%%%%%%%%%%%%%%%%%%%%%%%%%%%%
%%%%%%%%%%%%%%%%%%%%%%%%%%%%%%%%%%%%%%%%%%%%%%%%%%%%%%%%%%%%%%%

\end{document}